\documentclass[12pt]{amsart}%
\usepackage{amsmath}
\usepackage{amssymb}
\usepackage{amsfonts}
\usepackage{graphicx}%
\providecommand{\U}[1]{\protect\rule{.1in}{.1in}}
\newtheorem{theorem}{Theorem}
\theoremstyle{plain}

\newtheorem{corollary}{Corollary}

\newtheorem{lemma}{Lemma}

\newtheorem{proposition}{Proposition}
\newtheorem{remark}{Remark}

\numberwithin{equation}{section}
\begin{document}
\title{A note on delay-inverse systems, I}
\author{Nikica Ugle\v{s}i\'{c}}
\address{Sv. Ante 9, 23287 Veli R\aa t, Hrvatska (Croatia)}
\email{uglesic@pmfst.hr.}
\thanks{This paper is in final form and no version of it will be submitted for
publication elsewhere.}
\date{March 7, 2025}

\begin{abstract}
A generalization of an inverse system in a category was recently introduced,
as well as that of the corresponding pro-category These so called the
delay-inverse systems and delay-pro-category could potentially yield a new
theory of (delay-) inverse systems as well as a kind of coarser abstract shape
theory. However, we have proven that, whenever an indexing set has cardinality
$\aleph_{n}$, $n\in%
\mathbb{N}
_{0}$, the potential new theory reduces, in its essence (the classification
and invariants), to the ordinary one.

\end{abstract}
\maketitle

\section{\textbf{Introduction}\smallskip}

An interesting generalization of an inverse system $\boldsymbol{X}%
=(X_{a},p_{aa^{\prime}},A)$ in a category $\mathcal{C}$, called a
\textit{delay-inverse system}, was recently defined and studied by V.
Matijevi\'{c} and L. R. Rubin, [4, 5]. Thereby, \textit{delay} means that the
strict commutativity condition on bonding morphisms, $p_{aa^{\prime}%
}p_{a^{\prime}a^{\prime\prime}}=p_{aa^{\prime\prime}}$, $a\leq a^{\prime}\leq
a^{\prime\prime}$, is shifted to an $a^{\ast}$, $a\leq a^{\ast}\leq a^{\prime
}\leq a^{\prime\prime}$, that is far enough from $a$. The same idea is then
applied in order to obtain new morphisms and an equivalence relation on their
sets. The appropriate theory should rise by studying an associated
delay-pro-category of such systems and corresponding morphisms. The
constructed delay-pro-category significantly enlarges the ordinary
pro-category by the object class as well as by sets of morphisms, even between
the ordinary (commutative) inverse systems. Furthermore, there are objects
(for instance, some topological spaces) that give rise of associated
delay-inverse systems. However, if the cardinality $card(A)=\aleph_{n}$, $n\in%
\mathbb{N}
_{0}$, we show that this \textit{delay-inverse systems theory} can be reduced
to the standard (commutative) one. (The general continuum hypothesis is
assumed!) Therefore, although it brings a useful tool, in the case of such
delay-inverse systems, the new theory (possibly of coarser shapes in
$\mathcal{C}$) does not yield a coarser classification of objects neither new
invariants comparing to the ordinary theory (ordinary abstract shapes). Still
the delay-inverse systems technique could be a new useful tool in studying
\textquotedblleft complex\textquotedblright objects (especially, some classes
of locally bad topological spaces) and with them associated systems (limits,
resolutions, expansions, \v{C}ech systems) that consist of simpler objects.

\section{\textbf{Category} $Dpro$-$\mathcal{C}$ \smallskip}

Our category terminology is that of [1]. Let us briefly recall the needed
definitions from [5] used in the construction of the \textquotedblleft
pro\textquotedblright-category $Dpro$-$\mathcal{C}$ on the delay-inverse
systems in an arbitrary category $\mathcal{C}$. The generalization of an
inverse system $\boldsymbol{X}$, to a delay one, is achieved by relaxing the
strict commutativity condition as follows:

$(\forall a\in A)(\exists a^{\ast}\geq a)(\forall a^{\prime\prime}\geq
a^{\prime}\geq a^{\ast})$ $p_{aa^{\prime}}p_{a^{\prime}a^{\prime\prime}%
}=p_{aa^{\prime\prime}}$.

\noindent Such an $a^{\ast}$ is hereby often referred as a
\textit{\textquotedblleft commutation index\textquotedblright} for $a$. Given
a pair of delay-inverse systems $\boldsymbol{X}$ and $\boldsymbol{Y}%
=(Y_{b},q_{bb^{\prime}},B)$ in $\mathcal{C}$, the generalization of a morphism
$(f,f_{b})$ (of $inv$-$\mathcal{C}$) of $\boldsymbol{X}$ to $\boldsymbol{Y}$,
in order to achieve a delay-morphism, is subjected to the following weaker condition:

$(\forall b\in B)(\exists b_{\ast}\geq b)(\forall b^{\prime}\geq b_{\ast
})(\exists a\geq f(b),f(b^{\prime}))(\forall a^{\prime}\geq a)$

$q_{bb^{\prime}}f_{b^{\prime}}p_{f(b^{\prime})a^{\prime}}=f_{b}%
p_{f(b)a^{\prime}}$.

\noindent The identity delay-morphisms and the composition of those morphisms
are defined as in the category $inv$-$\mathcal{C}$. A certain new category,
denoted by $Dinv$-$\mathcal{C}$, in this way is obtained, such that
$inv$-$\mathcal{C}$ has become its subcategory.

\noindent Given a pair of delay-morphisms $(f,f_{b}),(f^{\prime},f_{b}%
^{\prime}):\boldsymbol{X}\rightarrow\boldsymbol{Y}$, an equivalence relation
on the set $Dinv$-$\mathcal{C}(\boldsymbol{X},\boldsymbol{Y})$ is defined by

$(\forall b\in B)(\exists a_{b}\geq f(b),f^{\prime}(b))(\forall a\geq a_{b})$

$f_{b^{\prime}}p_{f(b)a}=f_{b}^{\prime}p_{f^{\prime}(b)a}$.

\noindent In such a case, we write down $(f,f_{b})\overset{d}{\sim}(f^{\prime
},f_{b}^{\prime})$, and the corresponding equivalence class

$[(f,f_{b})]\in Dinv$-$\mathcal{C}(\boldsymbol{X},\boldsymbol{Y})/\overset
{d}{\sim}$

\noindent is denoted by $\boldsymbol{f}:\boldsymbol{X}\rightarrow
\boldsymbol{Y}$. Finally, since the equivalence relation $\overset{d}{\sim}$
is compatible with the composition of morphisms in $Dinv$-$\mathcal{C}$, there
is the corresponding quotient category

$Dinv$-$\mathcal{C}(\boldsymbol{X},\boldsymbol{Y})/\overset{d}{\sim}$ $\equiv
Dpro$-$\mathcal{C}$.

\noindent Clearly, $pro$-$\mathcal{C}$ is a subcategory of $Dpro$%
-$\mathcal{C}$. \newpage

\section{$Dpro$-$\mathcal{C}$ \textbf{versus} $pro$-$\mathcal{C}$}

In this section we show that the basic properties, analogous to those of
ordinary (commutative) inverse systems (see [3], Ch.I., Sec. 1.1, 1.2, 1.3)
hold true for delay-inverse systems as well. The statements concerning the
indexing sets only remain valid without any change. The first proposition is a
delay-analogue of [3], Theorem I.1.1,

\begin{proposition}
\label{P1}Let $\boldsymbol{X}=(X_{a},p_{aa^{\prime}},A)$ be a delay-inverse
system in a category $\mathcal{C}$, and let a subset $A^{\prime}\subseteq A$
be cofinal in $(A,\leq)$. Then $\boldsymbol{X}^{\prime}=(X_{a},p_{aa^{\prime}%
},A^{\prime})$ is also a delay-inverse system and the restriction morphism
$\boldsymbol{i}:\boldsymbol{X}\rightarrow\boldsymbol{X}^{\prime}$ is an
isomorphism in $Dpro$-$\mathcal{C}$.
\end{proposition}

\begin{proof}
Since $(A^{\prime},\leq)$ is cofinal in $(A,\leq)$, it is directed. and
$\boldsymbol{X}^{\prime}=(X_{a},p_{gla^{\prime}},A^{\prime})$ is indeed a
delay-inverse system in $\mathcal{C}$. Further, it is obvious that the
restriction morphism

$\boldsymbol{i}=[(i,i_{a}=1_{X_{a}})]:\boldsymbol{X}\rightarrow\boldsymbol{X}%
^{\prime}$, $\quad i:A^{\prime}\hookrightarrow A$,

\noindent belongs to $Dpro$-$\mathcal{C}(\boldsymbol{X},\boldsymbol{X}%
^{\prime})$. Since $A^{\prime}$ is cofinal in $A$, there is a function

$j:A\rightarrow A^{\prime}$, $a\mapsto j(a)\geq a^{\ast}$,

\noindent where $a^{\ast}\in A$ is a \textquotedblleft commutation
index\textquotedblright\ (in $\boldsymbol{X}$) for $a\in A^{\prime}\subseteq
A$. Put

$j_{a}=p_{aj(a)}:X_{j(a)}\rightarrow X_{a}$, $a\in A^{\prime}$.

\noindent One readily sees that

$\boldsymbol{j}=[(j,j_{a})]:\boldsymbol{X}^{\prime}\rightarrow\boldsymbol{X}$

\noindent is a morphism of $Dpro$-$\mathcal{C}(\boldsymbol{X}^{\prime
},\boldsymbol{X})$. The verification that $\boldsymbol{ji}=1_{\boldsymbol{X}}$
and $\boldsymbol{ij}=1_{\boldsymbol{X}^{\prime}}$ in $Dpro$-$\mathcal{C}$
hold, is a trivial routine.
\end{proof}

We are now proving a needed delay-analogue of [3], Theorem I.1.2.

\begin{proposition}
\label{P2}Every delay-inverse system $\boldsymbol{X}=(X_{a},p_{aa^{\prime}%
},A)$ in $\mathcal{C}$ admits an isomorphic (in $Dpro$-$\mathcal{C}$)
delay-inverse system $\boldsymbol{Y}=(Y_{b},q_{bb^{\prime}},B)$ in
$\mathcal{C}$ such that $B$ is ordered and cofinite having $card(B)\leq
card(A)$, and each term $Y_{b}$ in $\boldsymbol{Y}$ is actually a term $X_{a}$
of $\boldsymbol{X}$, while each bonding morphism $q_{bb^{\prime}}$ of
$\boldsymbol{Y}$ is a bonding morphism $p_{aa^{\prime}}$ of $\boldsymbol{X}$.
\end{proposition}

\begin{proof}
If $A$ is finite or if there exists $\max A$, then the claim is trivially
true. Suppose that $A$ is infinite and without the maximal element. By
Proposition 1 and by Remark I.1.1 of [3], we may assume that $(A,\leq)$ is
antisymmetric. We define a new indxeing set $B$ by means of so called
Marde\v{s}i\'{c} trick. The elements of $B$ are all finite subsets $b\subset
A$ having maximal elements (which are unique in all $b$). Since $A$ is
infinite, $card(B)=card(A)$ holds. $B$ is ordered by

$b\leq b^{\prime}\Leftrightarrow b\subseteq b^{\prime}$.

\noindent It is obvious that $(B,\leq)$ is cofinite, while an easy
verification shows that it is directed as well. Notice also that $b\leq
b^{\prime}$ in $(B,\leq)$ implies $\max b\leq\max b^{\prime}$ in $(A,\leq)$. Put

$Y_{b}=X_{\max b}$, $b\in B$,

$q_{bb^{\prime}}=p_{\max b\max b^{\prime}}:Y_{b^{\prime}}=X_{\max b^{\prime}%
}\rightarrow X_{\max b}=Y_{b}$, $b\leq b^{\prime}.$

\noindent Let us prove that $\boldsymbol{Y}=(Y_{b},q_{bb^{\prime}},B)$ is a
delay-inverse system in $\mathcal{C}$. Given a $b=\{a_{1},\cdots,a_{n}\}\in
B$, let $a_{i}=\max b$, where $i$ is unique, $1\leq i\leq n\in%
\mathbb{N}
$. Let $a^{\ast}\geq a_{i}$ be a \textquotedblleft commutation
index\textquotedblright\ for $a_{i}$. Then

$(\forall a^{\prime\prime}\geq a^{\prime}\geq a^{\ast})$ $p_{a_{i}a^{\prime}%
}p_{a^{\prime}a^{\prime\prime}}=p_{a_{i}a^{\prime\prime}}$.

\noindent Since $a^{\ast}\geq a_{i}$, it follows that

$b^{\ast}=b\cup\{a^{\ast}\}\in B$, $b^{\ast}\geq b$, $\max b^{\ast}=a^{\ast}$.

\noindent Let $b^{\prime\prime}\geq b^{\prime}\geq b^{\ast}$ in $(B,\leq)$. Then

$a^{\prime\prime}\equiv\max b^{\prime\prime}\geq a^{\prime}\equiv\max
b^{\prime}\geq a^{\ast}=\max b^{\ast}$,

\noindent and consequently,

$p_{a_{i}a^{\prime}}p_{a^{\prime}a^{\prime\prime}}=p_{a_{i}a^{\prime\prime}}$,

\noindent which means

$q_{bb^{\prime}}q_{b^{\prime}b^{\prime\prime}}=q_{bb^{\prime\prime}}$

\noindent and shows that $\boldsymbol{Y}$ is a delay-inverse system in
$\mathcal{C}$. It remains to prove that $\boldsymbol{X}$ and $\boldsymbol{Y}$
are mutually isomorphic in $Dpro$-$\mathcal{C}$. Let us define

$f:B\rightarrow A$, $f(b)=\max b$,

$f_{b}\equiv1_{X_{\max b}}:X_{f(b)}=X_{\max b}\rightarrow X_{\max b}=Y_{b}$,
$b\in B$.

\noindent It trivially follows that $(f,f_{b})\in Dinv$-$\mathcal{C}%
(\boldsymbol{X},\boldsymbol{Y}),$ and thus $\boldsymbol{f}=[(f,f_{b})]\in
Dpro$-$\mathcal{C}(\boldsymbol{X},\boldsymbol{Y})$. Conversely, define

$g:A\rightarrow B$, $g(a)=b\equiv\{a\}$,

$g_{a}=1_{X_{a}}:Y_{g(a)}=X_{a}\rightarrow X_{a}$, $a\in A$.

\noindent As before, it is obvious that $(g,g_{a})\in Dinv$-$\mathcal{C}%
(\boldsymbol{Y},\boldsymbol{X}),$ and thus $\boldsymbol{g}=[(g,g_{a})]\in
Dpro$-$\mathcal{C}(\boldsymbol{Y},\boldsymbol{X})$.

\noindent In order to prove that

$\boldsymbol{gf}=1_{\boldsymbol{X}}$, \quad i.e. $(g,g_{a})(f,f_{b}%
)\overset{d}{\sim}(1_{A},1_{X_{a}})$,

\noindent given an $a\in A$, for a needed $a_{\ast}\geq a$ one may take a
\textquotedblleft commutation index\textquotedblright\ $a^{\ast}$ for $a$ in
$\boldsymbol{X}$. Then the conclusion follows trivially. Similarly, in order
to prove that

$\boldsymbol{fg}=1_{\boldsymbol{Y}}$, \quad i.e. $(f,f_{b})(g,g_{b}%
)\overset{d}{\sim}(1_{B},1_{Y_{b}})$,

\noindent given a $b\in B$, for a needed $b_{\ast}\geq b$ one may take a
\textquotedblleft commutation index\textquotedblright\ $b^{\ast}$ for $b$ in
$\boldsymbol{Y}$. Again the conclusion follows trivially. This completes the
proof of the proposition.
\end{proof}

The third proposition assures that a countable delay-inverse system can be
reduced to a delay-inverse sequence. In the ordinary case, one can find this
claim in [2], Exercises, A.!., p. 229., while a proof may be that of Lemma 9
of [6].

\begin{proposition}
\label{P3}Let $\boldsymbol{X}=(X_{a},p_{aa^{\prime}},A)$ be a delay-inverse
system in a category $\mathcal{C}$ having $A$ countably infinite, i.e.,
$card(A)=\aleph_{0}$. Then there exists a strictly increasing sequence
$(a_{n})$ in $(A,\leq)$ such that

$\boldsymbol{X}^{\prime}=(X_{n}^{\prime}=X_{a_{n}},p_{nn^{\prime}}^{\prime
}=p_{a_{n}a_{n^{\prime}}},%
\mathbb{N}
)$

\noindent is a delay-inverse system (sequence) and $\boldsymbol{X}^{\prime}$
is isomorphic with $\boldsymbol{X}$ in $Dpro$--$\mathcal{C}$.
\end{proposition}

\begin{proof}
Since $A$ is countable, one may consider (by forgetting\ the ordering) that

$A=\{a_{1,}a_{2},\ldots a_{i},\ldots\}$, $i\in%
\mathbb{N}
$.

\noindent Let us denote $a_{1}^{\prime}\equiv a_{1}$, and put

$A^{(1)}=\{a\in A\mid a\leq a_{1}^{\prime}\}\subseteq A.$

\noindent If $A^{(1)}=A$, then $a_{1}^{\prime}$ is the maximal element of $A$.
Thus, to obtain a desired sequence one may put

$\boldsymbol{X}^{\prime}=(X_{i}^{\prime}=X_{a_{1}^{\prime}},p_{ii^{\prime}%
}^{\prime}=1,%
\mathbb{N}
)$.

\noindent Then $\boldsymbol{X}^{\prime}$ is a rudimentary object of
$Ob(pro$-$\mathcal{C})\subseteq Ob(Dpro$-$\mathcal{C})$ and, obviously,
$\boldsymbol{X}^{\prime}\cong\boldsymbol{X}$ in $Dpro$-$\mathcal{C}$.

\noindent Now, assuming that $A$ has no maximal element, let us choose an
$a_{2}^{\prime}\in A\setminus A^{(1)}$ such that $a_{2}^{\prime}>a_{1}%
^{\prime}$ and $a_{2}^{\prime}\geq a_{2}$, and put

$A^{(2)}=\{a\in A\setminus A^{(1)}\mid a\leq a_{2}^{\prime}\}\subseteq
A\setminus A^{(1)}.$

\noindent Suppose that $A^{(1)}$, $\cdots$, $A^{(i)}$ are constructed
inductively, by $i\in%
\mathbb{N}
$, in the same manner. Choose an

$a_{i+1}^{\prime}\in A\setminus(A^{(1)}\cup\cdots\cup A^{(i)})$

\noindent such that $a_{i+1}^{\prime}>a_{i}^{\prime}$ and $a_{i+1}^{\prime
}\geq a_{i+1}$, and put

$A^{(i+1)}=\{a\in A\setminus(A^{(1)}\cup\cdots\cup A^{(i)}\mid a\leq
a_{i+1}^{\prime}\}\subseteq A\setminus(A^{(1)}\cup\cdots\cup A^{(i)}).$

\noindent Notice that

$A=\sqcup_{i\in%
\mathbb{N}
}A^{(i)}$

\noindent(disjoint union) and the subset

$A^{\prime}\equiv\{a_{i}^{\prime}\mid i\in%
\mathbb{N}
\}\subseteq A$

\noindent is cofinal in $A$ (with respect to the given ordering). Clearly,
$a_{i}^{\prime}<a_{i^{\prime}}^{\prime}$ if and only if $i<i^{\prime}$.
Further, $a\leq a^{\prime}$ in $A$ if and only if either there exists an $i\in%
\mathbb{N}
$ such that $a\leq a^{\prime}$ in $A^{(i)}$, or there exists a pair
$i,i^{\prime}\in%
\mathbb{N}
$, such that $a\in A^{(i)}$, $a^{\prime}\in A^{(i^{\prime})}$, $a\leq
a^{\prime}$ and $i<i^{\prime}$. Let us define

$X_{i}^{\prime}=X_{a_{i}^{\prime}}$, $i\in%
\mathbb{N}
$;

$p_{ii^{\prime}}^{\prime}=p_{a_{i^{\prime}}^{\prime}a_{i}^{\prime}%
}:X_{i^{\prime}}^{\prime}\rightarrow X_{i}^{\prime}$ , $i\leq i^{\prime}$.

\noindent Since $(A^{\prime},\leq)$ is cofinal in $(A,\leq)$, it follows that
the subsystem $\boldsymbol{X}^{\prime}=(X_{i}^{\prime},p_{ii^{\prime}}%
^{\prime},%
\mathbb{N}
)$ of $\boldsymbol{X}$ is a delay-inverse sequence in $\mathcal{C}$. Moreover,
by Proposition 1, the restriction morphism $\boldsymbol{i}:\boldsymbol{X}%
\rightarrow\boldsymbol{X}^{\prime}$ is an isomorphism of $Dpro$-$\mathcal{C}$.
\end{proof}

\noindent. \newpage

\section{\textbf{Reduction to commutativity} \smallskip}

We are now ready to prove the main fact. Let us remind that the general
continuum hypothesis is assumed,and consequently, the axiom of choice holds
true as well as the statement that every set can be well ordered. In the first
part we shall consider delay-inverse systems having the indexing sets of
cardinality $\aleph_{n}$, $n\in%
\mathbb{N}
_{0}$, where $\aleph_{n+1}=2^{\aleph_{n}}$.

Our main theorem reads as follows.

\begin{theorem}
\label{T1}Every delay-inverse system $\boldsymbol{X}=(X_{a},p_{aa^{\prime}%
},A)$ in a category $\mathcal{C}$, having $card(A)=\aleph_{n}$ for some $n\in%
\mathbb{N}
_{0}$, is isomorphic in $Dpro$-$\mathcal{C}$ to a cofinite inverse system
$\boldsymbol{Y}=(Y_{b},q_{bb^{\prime}},B)$ in $\mathcal{C}$. Moreover, every
$Y_{b}$ is some $X_{a}$ and every $q_{bb^{\prime}}$ is some $p_{aa^{\prime}}$.
\end{theorem}

At first we consider the simplest case of an $A$ having cardinality
$card(A)\equiv$ $\left\vert A\right\vert \leq\aleph_{0}$, and immediately
after that the case of $\left\vert A\right\vert =2^{\aleph_{0}}$.The prof of a
countable case reduces to a proof for a sequence. Our proof runs
straightforwardly by iteration of the \textquotedblleft commutation
indices\textquotedblright. The proof of the smallest uncountable case combines
simultaneously the ordinary (by $%
\mathbb{N}
_{0}$) and transfinite (by $\left[  0,\omega_{1}\right\rangle $) induction.

\begin{lemma}
\label{L1}Every delay-inverse system $\boldsymbol{X}=(X_{a},p_{aa^{\prime}%
},A)$ in a category $\mathcal{C}$, satisfying $card(A)\leq2^{\aleph_{0}}$, is
isomorphic in $Dpro$-$\mathcal{C}$ to a cofinite inverse system
$\boldsymbol{Y}=(Y_{b},q_{bb^{\prime}},B),$in $\mathcal{C}$. Moreover, every
$Y_{b}$ is some $X_{a}$ and every $q_{bb^{\prime}}$ is some $p_{aa^{\prime}}$.
\end{lemma}

\begin{proof}
We split our proof into two parts according to cardinality of $A$.

\textbf{1. The countable case.}

\noindent If the indexing set $A$ is finite or $(A,\leq)$ (generally) has the
maximal element, the proof is trivial. Let $A$ be countably infinite, i.e.,
$\left\vert A\right\vert =\aleph_{0}$ (and without maximum). By Proposition 3,
we may assume that $A=%
\mathbb{N}
$. So let $\boldsymbol{X}=(X_{j},p_{jj^{\prime}},%
\mathbb{N}
)$ be a delay-inverse sequence in $\mathcal{C}$. Take $j\equiv j_{1}=1$, and
let $j_{1}^{\ast}\geq j_{1}$ be the minimal \textquotedblleft commutation
index\textquotedblright\ for $j_{1}$ in $\boldsymbol{X}$. Then

$(\forall j^{\prime\prime}\geq j^{\prime}\geq j_{1}^{\ast})$ $p_{j_{1}%
j^{\prime}}p_{j^{\prime}j^{\prime\prime}}=p_{j_{1}j^{\prime\prime}}.$

\noindent Put $j_{2}=j_{1}^{\ast}+1$. Then $j_{2}>j_{1}$, and clearly again

$(\forall j^{\prime\prime}\geq j^{\prime}\geq j_{2})$ $p_{j_{1}j^{\prime}%
}p_{j^{\prime}j^{\prime\prime}}=p_{j_{1}j^{\prime\prime}}$

\noindent and especially

$p_{j_{1}j_{2}}p_{j_{2}j^{\prime\prime}}=p_{j_{1}j^{\prime\prime}}$.

\noindent hold. In the same way, let $j_{2}^{\ast}\geq j_{2}$ be the minimal
\textquotedblleft commutation\textquotedblright\ index for $j_{2}$ in
$\boldsymbol{X}$. Put $j_{3}=j_{2}^{\ast}+1$. Then

$(\forall j^{\prime\prime}\geq j^{\prime}\geq j_{3})$ $p_{j_{2}j^{\prime}%
}p_{j^{\prime}j^{\prime\prime}}=p_{j_{2}j^{\prime\prime}}$

\noindent holds, and especially

$p_{j_{2}j_{3}}p_{j_{3}j^{\prime\prime}}=p_{j_{2}j^{\prime\prime}}$.

\noindent According to the previous step,

$p_{j_{1}j_{2}}p_{j_{2}j_{3}}=p_{j_{1}j_{3}}$.

\noindent holds as well. Continuing inductively, let $j_{3}^{\ast}\geq j_{3}$
be the minimal \textquotedblleft commutation index\textquotedblright\ for
$j_{3}$ in $\boldsymbol{X}$, and put $j_{4}=j_{3}^{\ast}+1$. By this inductive
iteration one obtains a strictly increasing sequence $(j_{k})$ in $%
\mathbb{N}
$, $j_{1}=1$, such that

$(\forall j_{k}\leq j_{k^{\prime}}\leq j_{k^{\prime\prime}})$ $p_{j_{k}%
j_{k^{\prime}}}p_{j_{k^{\prime}}j_{k^{\prime\prime}}}=p_{j_{k}j_{k^{\prime
\prime}}}$.

\noindent Notice that $j_{k}\leq j_{k^{\prime}}\leq j_{k^{\prime\prime}}$ if
and only if $k\leq k^{\prime}\leq k^{\prime\prime}$. Put

$X_{k}^{\prime}=X_{j_{k}}$, $p_{kk^{\prime}}^{\prime}=p_{j_{k}j_{k^{\prime}}}%
$, $k\in%
\mathbb{N}
$.

\noindent In this way an inverse sequence $\boldsymbol{X}^{\prime}%
=(X_{k}^{\prime},p_{kk^{\prime}}^{\prime},%
\mathbb{N}
)$ (commutative) is constructed, i.e. $\boldsymbol{X}^{\prime}\in
Ob(pro$-$\mathcal{C})$, which is a strictly increasing subsequence of
$\boldsymbol{X}\in Ob(Dpro$-$\mathcal{C})$. By Proposition 1, the restriction
morphism $\boldsymbol{i}:\boldsymbol{X}\rightarrow\boldsymbol{X}^{\prime}$ is
an isomorphism in $Dpro$-$\mathcal{C}$. This completes our proof in the
countable case.

\textbf{2. The }$\aleph_{1}=2^{\aleph_{0}}$\textbf{-uncountable case.}

Let $A$ be uncountable having $\left\vert A\right\vert =2^{\aleph_{0}}$. By
Proposition 2 (\textquotedblleft Marde\v{s}i\'{c} trick\textquotedblright),
one obtains the cofinite delay-inverse system of the same \textquotedblleft
material\textquotedblright\ and with an indexing set of the same cardinality
$2^{\aleph_{0}}$, that is isomorphic in $Dpro$-$\mathcal{C}$ to
$\boldsymbol{X}$. Thus we may assume, without loss of generality, that
$\boldsymbol{X}$ is already such a system. Let us denote

$A=\{a_{j}\mid j\in J\}$, $\left\vert J\right\vert =2^{\aleph_{0}}$,

\noindent and further

$J=\sqcup_{k\in%
\mathbb{N}
_{0}}J_{k}$, $\quad J_{k}=\{j\mid\left\vert a_{j}\right\vert =k\}\subseteq J$;

$A=\sqcup_{k\in%
\mathbb{N}
_{0}}A_{k}$, $\quad A_{k}=\{a_{j}\mid j\in J_{k}\}\subseteq A$.

\noindent Hereby $\left\vert a_{j}\right\vert $ denotes the number of all
$a\in A$ such that $a<a_{j}$ (i.e., $a\leq a_{j}$ and $a\neq a_{j}$). Since
the set $A$ is directed and infinite, for every $k\in%
\mathbb{N}
_{0}$,

$\left\vert A_{k}\right\vert =$ $\left\vert J_{k}\right\vert =$ $\left\vert
J\right\vert =$ $\left\vert A\right\vert \equiv2^{\aleph_{0}}$.

\noindent Notice that, for each $k\in%
\mathbb{N}
_{0}$, any pair $a,a^{\prime}\in A_{k}$, $a\neq a^{\prime}$, is not related in
$(A,\leq)$. Let $\leq_{k}$ be a well order on $J_{k}$, $k\in%
\mathbb{N}
_{0}$, and let us write $(J_{k},\leq_{k})$ down by following the order as follows:

$j_{0}^{k},j_{1}^{k},\cdots,j_{n}^{k},\cdots,j_{\omega}^{k},j_{\omega+1}%
^{k},\cdots,j_{\omega+\omega}^{k},\cdots,\cdots,j_{s}^{k},j_{s+1}^{k},\cdots$

\noindent where $\omega$ is the first infinite ordinal number. Since
$\left\vert J_{k}\right\vert =2^{\aleph_{0}}$, one may assume with no loss of
generality that every $(J_{k},\leq_{k})$ is isomorphic to ($\left[
0.\omega_{1}\right\rangle ,\leq)$, where $\omega_{1}$ is the first uncountable
ordinal. We remind that there is no relation between the ordering $\leq$ on
$A$ to any well ordering $\leq_{k}$ on $J_{k}$. Recall also that, for each
limit ordinal $r$, $\omega\leq r<\omega_{1}$, there exists the
\textquotedblleft next\textquotedblright\ (closest greater) limit ordinal,
which is $r+\omega$. Our proof proceeds by combining the transfinite induction
over the set of all limit ordinals of $\left[  0,\omega_{1}\right\rangle $
followed, in each step, by the ordinary (over $%
\mathbb{N}
_{0}$) inductions.

\noindent We begin with $(J_{0},\leq_{0})$ by considering the subset $A_{0}$
consisting of all the elements $a_{j_{s}^{9}}\equiv a_{j_{s}^{o}}^{0}\in A$
which have no predecessors, i.e., $\left\vert a_{j_{s}^{0}}^{0}\right\vert
=0$. Consider the first element $j_{9}^{0}\in J_{0}$, i.e., the corresponding
index $a_{j_{0}^{9}}^{0}\in A_{0}\equiv A_{0}^{0}$. Since $\boldsymbol{X}$ is
a delay-inverse system, there exists an $a^{\ast}\in A$, $a^{\ast}\geq
a_{j_{0}^{9}}^{0}$, such that, for all $a^{\prime\prime}\geq a^{\prime}\geq
a^{\ast}$, the condition

$p_{a_{j_{0}^{0}}^{0}a^{\prime}}p_{a^{\prime}a^{\prime\prime}}=p_{a_{j_{0}%
^{0}}^{0}a^{\prime\prime}}$

\noindent is fulfilled. Choose one of such $a^{\ast}$ (for instance, of
minimal $\left\vert a^{\ast}\right\vert \geq0$), and then, by directedness of
$(A,\leq),$ there exists an $\overline{a}\in A$, such that $\overline{a}\geq
a^{\ast},a_{j_{0}^{1}}$ ($a_{j_{0}^{1}}\in A_{1}$ - the first element by
$\leq_{1}$ of $J_{1}$). Pick up one of such $\overline{a}$, for instance, of
minimal $\left\vert \overline{a}\right\vert \geq1$, and denote $a_{j_{0}^{0}%
}^{1}\equiv\overline{a}$.

\noindent In the same way, for $a_{j_{1}^{0}}^{0}\in A_{0}$, there exists its
\textquotedblleft commutation index\textquotedblright\ $a^{\ast}\in A$,
$a^{\ast}\geq a_{j_{1}^{0}}^{0}$, and further, there exists an $\overline
{a}\in A$, $\overline{a}\geq a^{\ast},a_{j_{0}^{0}}^{1},a_{j_{1}^{1}}$. Choose
such an $\overline{a}$ (for instance, of minimal $\left\vert \overline
{a}\right\vert $) and denote $a_{j_{1}^{0}}^{1}\equiv\overline{a}$.

\noindent Assume, by induction on $n\in%
\mathbb{N}
_{0}$ ( in this first step of the transfinite induction), that, given an $n$,
the indices $a_{j_{i}^{0}}^{1}$, $0\leq i\leq n$, are chosen, satisfying
$\left\vert a_{j_{i}^{0}}^{1}\right\vert \geq1$, such that, for every pair
$a_{j_{i}^{0}}^{0},a_{j_{i^{\prime}}^{0}}^{0}$, $0\leq i,i^{\prime}\leq n$,
there is an $i^{\prime\prime}\leq n$, such that $a_{j_{i^{\prime\prime}}^{0}%
}^{1}\geq a_{j_{i}^{0}}^{0},a_{j_{i^{\prime}}^{0}}^{0}$, and that, for every
$a_{j_{i}^{1}}$, $0\leq i\leq n$, $a_{j_{i}^{0}}^{1}\geq a_{j_{i}^{1}}$.
(Hereby, $i^{\prime\prime}=\max\{i,i^{\prime}\}$.) Now, as in the step
$0\longmapsto1$, for $a_{j_{n+1}^{0}}^{0}$ choose firstly a \textquotedblleft
commutation index\textquotedblright\ $a^{\ast}\geq a_{j_{n+1}^{0}}^{0}$, and
then an $\overline{a}\geq a^{\ast},a_{j_{n}^{0}}^{1},a_{j_{n+1}^{1}}$, and
denote $a_{j^{0}}^{1}\equiv\overline{a}$. In this way the mentioned properties
of $a_{j_{i}^{0}}^{1}$, $0\leq i\leq n$, are prolonged to $i=n+1$. The
induction by $n\in%
\mathbb{N}
_{0}$ yields a sequence $(a_{j_{n}^{0}}^{1})$ (increasing in $(\sqcup_{k\geq
1}A_{k},\leq)$ $\subseteq(A,\leq)$) having, with respect to sequences
$(a_{j_{n}^{0}}^{0})$ and $(a_{j_{n}^{1}})$, the following properties:

\noindent(i)$\quad(\forall n,n^{\prime}\in%
\mathbb{N}
_{0})(\exists n^{\prime\prime}\in%
\mathbb{N}
_{0})$ $a_{j_{n^{\prime\prime}}^{0}}^{1}\geq a_{j_{n}^{0}}^{0},a_{j_{n^{\prime
}}^{0}}^{0}$;

\noindent(ii)\quad$(\forall n\in%
\mathbb{N}
_{0})(\exists n^{\prime}\in%
\mathbb{N}
_{0})$ $a_{j_{n^{\prime}}^{0}}^{1}\geq a_{j_{n}^{1}}$;

\noindent(iii)\quad$(\forall a\leq a^{\prime}\leq a^{\prime\prime}\in
A_{0}^{0}\cup\{a_{j_{n}^{0}}^{1}\mid n\in%
\mathbb{N}
_{0}\})$

$p_{aa^{\prime}}p_{a^{\prime}a^{\prime\prime}}=p_{aa^{\prime\prime}}$.

\noindent Notice that, at this starting level of the first step, in (i) it is
$n^{\prime\prime}=\max\{n,n^{\prime}\}$, while in (ii), $n^{\prime}=n$.

\noindent Consider now the first (minimal) limit ordinal number $\omega$, i.e.
$j_{\omega}^{0}$ of $(J_{o},\leq_{0})$, and the corresponding index
$a_{j_{\omega}^{0}}^{0}\in A_{0}$. Choose a \textquotedblleft commutation
index\textquotedblright\ $a^{\ast}\geq a_{j_{\omega}^{0}}^{0}$, and then an
$\overline{a}\geq a^{\ast},a_{j_{\omega}^{1}},a_{j_{0}^{0}}^{1}$. Clearly,
$\left\vert \overline{a}\right\vert \geq1$. Denote $a_{j_{\omega}^{0}}%
^{1}\equiv\overline{a}$. Further, for $a_{j_{\omega+1}^{0}}^{0}$, choose its
$a^{\ast}\geq a_{j_{\omega+1}^{0}}^{0}$, and then an $\overline{a}\geq
a^{\ast},a_{j_{\omega+1}^{1}},a_{j_{\omega}^{0}}^{1},a_{j1}^{1}$, and denote
$\overline{a}\equiv a_{j_{\omega+1}^{0}}^{1}$. By mimicking the previous
inductive construction (over $%
\mathbb{N}
_{0}$, and carrying about the sequence $(a_{j_{n}^{0}}^{1})$), one obtains a
sequence $(a_{j_{\omega+n}^{0}}^{1})$ (increasing in $(\sqcup_{k\geq1}%
A_{k},\leq)$ $\subseteq(A,\leq)$) having properties similar to (i), (ii) and
(iii), now with respect to $(a_{j_{\omega+n}^{0}}^{0})$ and $(a_{j_{\omega
+n}^{1}})$. Moreover, they all together, i.e. $(a_{j_{s}^{0}}^{1})$, $0\leq
s<\omega+\omega$, have those properties (necessary for the needed
directedness, cofinality and commutativity) with respect to $(a_{j_{s}^{0}%
}^{0})$ and $(a_{j_{s}^{1}})$, $0\leq s<\omega+\omega$.

\noindent Let us now make the first transfinite inductive step from $\omega$
to the second limit ordinal $\omega+\omega$. Consider the index $j_{\omega
+\omega}^{0}\in J_{0}$ and the corresponding $a_{j_{\omega+\omega}^{0}}^{0}\in
A_{0}^{0}\equiv A_{0}$. Firstly, choose a \textquotedblleft commutation
index\textquotedblright\ $a^{\ast}$ for $a_{j_{\omega+\omega}^{0}}^{0}$, and
then an $\overline{a}\geq a^{\ast},a_{j_{\omega+\omega}^{1}},a_{j_{\omega}%
^{0}}^{1}$. Denote $a_{j_{\omega+\omega}^{0}}^{1}\equiv\overline{a}$ and then
proceed, analogously to the previous case of $\omega$, with the inductive
construction (over $%
\mathbb{N}
_{0}$) of an appropriate increasing sequence $(a_{j_{\omega+\omega+n}^{0}}%
^{1})$, $n\in%
\mathbb{N}
_{0}$.

\noindent Keeping in mind that our construction is based on intervals $\left[
r,r+\omega\right\rangle $, where $r<\omega_{1}$ is a limit ordinal, we can now
proceed by transfinite induction. Let $r$, $\omega+\omega\leq r<\omega_{1}$,
be an arbitrary limit ordinal, and assume that the indices $a_{j_{s}^{0}}%
^{1}\in A$, $\left\vert a_{j_{s}^{0}}^{1}\right\vert $ $\geq1$, are chosen for
all $s<r$, such that they fulfill conditions corresponding to (i), (ii) and
(iii) with respect to $(a_{j_{s}^{0}}^{0}$ and $(a_{j_{s}^{1}})$. Since the
interval $\left[  0,r\right\rangle $ is countable, its subset $R$ of all limit
ordinals $s$ less than $r$ is also countable. We may assume, without loss of
generality, that $R$, as a sequence $(s_{n})$, is not a stationary one, even
more, that all of its members are mutually different. However, clearly, this
sequence (a bijection of $%
\mathbb{N}
_{0}$ onto $R$)

$s_{0}=\omega,s_{1},\cdots,s_{n},\cdots<r$

\noindent cannot, generally, retains the canonical well order of $\left[
0,\omega_{1}\right\rangle $. In the same way, the subsets of all $a_{j_{s}%
^{0}}^{0}$, all $a_{j_{s}^{1}}$ and all $a_{j_{s}^{0}}^{1}$, $s\in\{s_{n}\mid
n\in%
\mathbb{N}
_{0}\}$, admit the sequence presentations (without the canonical well order)
as follows:

$a_{j_{s_{0}}^{0}}^{0}=a_{j_{0}^{0}}^{0},a_{j_{s_{1}}^{0}}^{0},\cdots
,a_{j_{s_{n}}^{0}}^{0},\cdots;$

$a_{j_{s_{0}}^{1}}=a_{j_{0}^{1}},a_{j_{s_{1}}^{1}},\cdots,a_{j_{s_{n}}^{1}%
},\cdots;$

$a_{j_{s_{0}}^{0}}^{1}=a_{j_{0}^{0}}^{1},a_{j_{s_{1}}^{0}}^{1},\cdots
,a_{j_{s_{n}}^{0}}^{1},\cdots.$

\noindent Let us repeat the technique used in the case $\omega+\omega$,
applied firstly to $r$ and to the first members of the three sequences from
above, and further (by induction on $%
\mathbb{N}
_{0}$) over the interval $\left[  r,r+\omega\right\rangle $. It yields a set

$A_{0r}^{1}\equiv\{a_{j_{s}^{0}}^{1}\mid j_{s}^{0}\in J_{0}$, $s<r+\omega
\}\subseteq\sqcup_{k\geq1}A_{k}\subset A$

\noindent such that the following conditions:

\noindent(I)$\quad(\forall s,s^{\prime}<r+\omega)(\forall a_{j_{s}^{0}}%
^{0},a_{j_{s^{\prime}}^{0}}^{0}\in A_{0}^{0})(\exists s^{\prime\prime
}<r+\omega$, i.., $a_{j_{s^{\prime\prime}}^{0}}^{1}\in A_{0r}^{1})$

$a_{j_{s^{\prime\prime}}^{0}}^{1}\geq a_{j_{s}^{0}}^{0},a_{j_{s^{\prime}}^{0}%
}^{0}$;

\noindent(II)\quad$(\forall s<r+\omega$, i.e., $a_{j_{s}^{1}}\in
A_{1})(\exists s^{\prime}<r+\omega$, i.e., $a_{j_{s^{\prime}}^{0}}^{1}\in
A_{0r}^{1})$

$a_{j_{s^{\prime}}^{0}}^{1}\geq a_{j_{s}^{1}}$;

\noindent(III)\quad$(\forall a\leq a^{\prime}\leq a^{\prime\prime}\in
A_{0}^{0}\cup A_{0r}^{1})$

$p_{aa^{\prime}}p_{a^{\prime}a^{\prime\prime}}=p_{aa^{\prime\prime}}$.

\noindent are fulfilled. Namely in order to achieve all needed properties of
chosen indices (at this level), one has to pass throughout the whole interval
$\left[  r,r+\omega\right\rangle $. Then, for instance (for directedness),
given a pair $a_{j_{s}^{0}}^{0},a_{j_{s^{\prime}}^{0}}^{0}$, where
$s\in\left[  s_{n},s_{n}+\omega\right\rangle $ and $s^{\prime}\in\left[
s_{^{\prime}n},s_{n^{\prime}}+\omega\right\rangle $, there exists an
$s_{n^{\prime\prime}}\geq s_{n},s_{n^{\prime}}$ and, by the construction,
there is an $s^{\prime\prime}\in\left[  s_{n^{\prime\prime}},s_{n^{\prime
\prime}}+\omega\right\rangle $ such that the chosen index $a_{j_{s^{\prime
\prime}}^{0}}^{1}$ satisfies (I) (see also Remark 1 below). Further, condition
(II) (for cofinality) is fulfilled because we always choose an $a^{\ast}\geq
a_{j_{s}^{0}}^{0}$, $a_{j_{s}^{0}}^{0}\in A$, to be greater than the
corresponding index $a_{j_{s}^{1}}\in A_{1}$ as well. Finally, condition (III)
(for commutativity) is fulfilled by choosing a \textquotedblleft commutation
index\textquotedblright\ of an $a_{j_{s}^{0}}^{0}$ to be valid for all
relevant predecessors $a_{1}^{\prime},\cdots,a_{m}^{\prime}$ of $a_{j_{s}^{0}%
}^{0}$ as well.

\noindent Thus, at this step of our proof, the set of indices

$A_{0}^{1}=\cup_{r}A_{0r}^{1}\subseteq\sqcup_{k\geq1}A_{k}$,$\subset A$

\noindent is constructed (by transfinite induction) such that it fulfills
conditions (I), (II) and (III) with respect to whole $A_{0}^{0}$ and $A_{1}$.
This completes the proof of the basic step, for all $j_{s}\in J_{0}$,
$s<\omega_{1}$. We can now push our construction from $A_{0}\equiv A_{0}^{0}$
to $A_{0}^{1}$.

\noindent By mimicking, and carrying about one additional condition, the whole
construction in the first step (where we were dealing with $A_{0}\equiv
A_{0}^{0}$ and $A_{1}$ in order to obtain $A_{0}^{1}$), applied now to $A_{1}$
and $A_{0}^{1}$, one obtains a set

$A_{1}^{2}\subseteq\sqcup_{k\geq2}A_{k}\subset A$

\noindent satisfying the analogues of (I), (II) and (III) with respect to
$A_{0}^{1}$ and $A_{2}$. The additional condition is as follows:

\noindent Let $a_{j_{s}^{0}}^{1}\in A_{0}^{1}$. Consider the set (finite,
because $A$ is cofinite) of all of its predecessors, and among them those
which are not related to $a_{j_{s}^{0}}^{0}$, i.e., those which do not belong
to any segment $[a_{j_{s}^{0}}^{0},a_{j_{s}^{0}}^{1}]$ of $(A,\leq)$. (Namely,
they have become out of the \textquotedblleft game\textquotedblright!) Denote
the considering predecessors by $a_{1}^{\prime},\cdots a_{m}^{\prime}$. Now,
when one chooses a \textquotedblleft commutation index\textquotedblright%
\ $a^{\ast}$ for $a_{j_{s}^{0}}^{1}$, let it be enough large such it can be a
\textquotedblleft commutation index\textquotedblright\ for every
$a_{i}^{\prime}$, $1\leq i\leq m$, as well. Then $a^{\ast}$ is also a
\textquotedblleft commutation index\textquotedblright\ for each relevant
predecessor of $a_{j_{s}^{0}}^{1}$.

\noindent Continuing inductively by $k\in%
\mathbb{N}
_{0}$, given an arbitrary $k_{0}\geq0$, suppose that the sets

$A_{m}^{m+1}\subseteq\sqcup_{k\geq m+1}A_{k}$, $0\leq m\leq k_{0}$,

\noindent are constructed such that the appropriate analogues of (I), (II) and
(III) hold also true \textquotedblleft gradually\textquotedblright. Now, in
the manner of the first step and with the additional condition, we repeat the
transfinite inductive construction hereby applied to $A_{k_{0}}^{k_{0}+1}$ and
$A_{k_{0}+1}$. Then the resulting subset

$A_{k_{0}+1}^{k_{0}+2}\subseteq\sqcup_{k\geq k_{0}+2}A_{k}\subset A$

\noindent satisfies the analogues of (I), (II) and (III) with respect to
$A_{k_{0}}^{k_{0}+1}$ and $A_{k_{0}+1}$ as well. With this, the whole
inductive construction is finished. A desired cofinite inverse system
(commutative) may be written as

$\boldsymbol{Y}=(Y_{b},q_{\mu\mu^{\prime}},B)$, $\quad B=A_{0}^{0}\cup
(\cup_{k\in%
\mathbb{N}
_{0}}A_{k}^{k+1})$,

\noindent where each $Y_{b}$ is the corresponding $X_{a}$, while each
$q_{\mu\mu^{\prime}}$ is the corresponding $p_{aa^{\prime}}$. Namely, the
directedness of $B$ holds by (I), the cofinality of $B$ in $A$ - by (II), the
cofinitness of $B$ follows by that of $A$, while the commutativity of
$\boldsymbol{Y}$ holds by choosing a \textquotedblleft commutation
index\textquotedblright\ of an $a$ to be greater than \textquotedblleft
commutation indices\textquotedblright\ of all (finitely many) relevant
predecessors of $a$ in $A$. Then, finally, $\boldsymbol{Y}\cong\boldsymbol{X}$
in $Dpro$.$\mathcal{C}$ follows by Proposition 1.
\end{proof}

\begin{remark}
\label{R1}Notice that the proof of Lemma 1 (the countable case) covers also
Proposition 3. Further, in the proof of the uncountable case, it was
\textquotedblleft insisted\textquotedblright\ on the directedness, conditions
(i) and (I), of the set of all before chosen indices. In order to achieve
this, we have used the fact that, for every ordinal $r\in\left[  0,\omega
_{1}\right\rangle $, there are at most countable many predecessors $s<r$.
However, conditions (i) and (I) may be abandoned throughout the proof, Namely,
in each step, the directedness in the \textquotedblleft
block\textquotedblright\ (an increasing sequence) is trivially fulfilled,
while the needed (only) final directedness follow by the cofinality of $B$ in
$A$ and the directedness of $A$.
\end{remark}

\begin{proof}
(of Theorem 1) Let the cardinality

$\left\vert A\right\vert \in\{\aleph_{0},\aleph_{n+1}=2^{\aleph_{n}}\mid n\in%
\mathbb{N}
_{0}\}$.

\noindent If $\left\vert A\right\vert =\aleph_{0}$, the statement follows by
the countable case of Lemma 1. Let

$\left\vert A\right\vert \in\{\aleph_{n+1}=2^{\aleph_{n}}\mid n\in%
\mathbb{N}
_{0}\}$.

\noindent The proof is by ordinary induction and it is based on applying Lemma
1 (and its proof) to each considered case, i.e., in each step. First of all,
Lemma 1 provides the proof of the basic step $\left\vert A\right\vert
=\aleph_{1}=2^{\aleph_{0}}$ , i.e., if $n=0$, there is an inverse system
$\boldsymbol{Y}$ that is a subsystem of the delay-inverse system
$\boldsymbol{X}$ such that $\boldsymbol{Y}\cong\boldsymbol{X}$ in
$Dpro$-$\mathcal{C}$. For our purpose, denote it by

$\boldsymbol{Y}_{0}=(Y_{b},q_{bb^{\prime}},B_{0})$,

\noindent where each $Y_{b}$ is the corresponding $X_{a}$ and each
$q_{bb^{\prime}}$ is the corresponding $p_{aa^{\prime}}$. Let $n=1$. i.e., let
$\boldsymbol{X}$ has an $A$ such that $\left\vert A\right\vert =\aleph
_{2}=2^{\aleph_{1}}$. By the basic step, the restriction of $\left[
0,\omega_{2}\right\rangle $ to the initial interval $\left[  0,\omega
_{1}\right\rangle $ yields (as in the proof of Lemma 1) a subset of indices

$B_{0}^{1}=A_{0}^{01}\cup(\cup_{k\in%
\mathbb{N}
_{0}}A_{k}^{k+1,1})\subseteq A$

\noindent such that the collection

$\boldsymbol{Y}_{0}^{1}=(Y_{b},q_{bb^{\prime}},B_{0}^{1})$

\noindent is an inverse system (subsystem of $\boldsymbol{X}$) having its
initial indices subordinated to the elements of the initial interval $\left[
0,\omega_{1}\right\rangle $ of $\left[  0,\omega_{2}\right\rangle $. Do the
same on the next interval $\left[  \omega_{1},\omega_{1}+\omega_{1}%
\right\rangle \subseteq\left[  0,\omega_{2}\right\rangle $, and go on. Notice
that $\left[  0,\omega_{2}\right\rangle $ consists of $\aleph_{2}$-many
consecutive intervals $\left[  r_{s},r_{s}+\omega_{1}\right\rangle $
isomorphic to $\left[  0,\omega_{1}\right\rangle $, where $r_{0}=0$,
$r_{1}=\omega_{1}$ is the limit ordinal of $\left[  0,\omega_{1}\right\rangle
$ in $\left[  0,\omega_{2}\right\rangle $ and $r_{s+1}$ is the limit ordinal
of $\left[  r_{s},r_{s}+\omega_{1}\right\rangle $ in $\left[  0,\omega
_{2}\right\rangle $. Thus, for each $s$, by Lemma 1 (see also Remark 1)
applied to each $\left[  r_{s},r_{s}+\omega_{1}\right\rangle $, there exists
an ordered subset $B_{s}^{1}\subseteq A$ which is cofinite and cofinal in the
subset of $A$ having the initial indices subordinated to elements of $\left[
r_{s},r_{s}+\omega_{1}\right\rangle $. (This does not involve any induction
because one uses only the fact that every ordinal has the immediate successor,
i.e., the construction for an $s$ does not depend on constructions in the
proceeding steps!) Denote by

$\boldsymbol{Y}_{s}^{1}=(Y_{b},q_{bb^{\prime}},B_{s}^{1})$

\noindent the corresponding collection of terms of $\boldsymbol{X}$. Then
$\boldsymbol{Y}_{s}^{1}$ is an inverse subsystem of $\boldsymbol{X}$, which is
commutative by the construction made in the proof of Lemma 1. Put

$\boldsymbol{Y}_{1}=(Y_{b},q_{\mu\mu^{\prime}},B_{1})$, $\quad B_{1}=B_{0}%
^{1}\cup(\cup B_{s}^{1})\subseteq A$,

\noindent where each $Y_{b}$ is the corresponding $X_{a}$ and each
$q_{bb^{\prime}}$ is the corresponding $p_{aa^{\prime}}$. Then $\boldsymbol{Y}%
_{1}$ is a desired inverse subsystem of $\boldsymbol{X}$ in the case of $n=1$.
Indeed, firstly, the subset $B_{1}$ is ordered, cofinite and cofinal in $A$.
Furthermore, it is directed as well. Namely, if $b,b^{\prime}\in B_{s}^{1}$,
for some $s$, then Lemma 1 provides a $b^{\prime\prime}\in B_{s}^{1}$,
$b^{\prime\prime}\geq b,b^{\prime}$. If $b\in B_{s}^{1}$ and $b^{\prime}\in
B_{s^{\prime}}^{1}$, $s\neq s^{\prime}$, then, since $A$ is directed, there is
an $a\in A$, $a\geq b,b^{\prime}$.Since $B_{1}$ is cofinal in $A$, there
exists a $b^{\prime\prime}\in B_{1}$, $b^{\prime\prime}\geq a$, and
consequently, $b^{\prime\prime}\geq b,b^{\prime}$. Finally, $\boldsymbol{Y}%
_{1}$ is commutative because of the way of choosing \textquotedblleft
commutation indices\textquotedblright\ (for all relevant predecessors as well)
of the considered elements of $A$\ in $\boldsymbol{X}$.

\noindent Now, given an $n\in%
\mathbb{N}
_{0}$, suppose that the statement of the theorem holds true for every $i$,
$0\leq i\leq n$. Then one can, analogously to the passing $0\mapsto1$, obtain
an inverse system $\boldsymbol{Y}_{n+1}$ by \textquotedblleft
gluing\textquotedblright\ $\aleph_{n+1}=2^{\aleph_{n}}$-many \textquotedblleft
copies\textquotedblright\ of $\boldsymbol{Y}_{n}$ throughout the all intervals

$\left[  r_{s},r_{s}+\omega_{n}\right\rangle \subseteq$ $\left[
0,\omega_{n+1}\right\rangle $,

\noindent where $r_{0}=0$, $r_{1}=\omega_{n}$ is the limit ordinal of $\left[
0,\omega_{n}\right\rangle $ in $\left[  0,\omega_{n+1}\right\rangle $ and
$r_{s+1}$ is the limit ordinal of $\left[  r_{s},r_{s}+\omega_{n}\right\rangle
$ in $\left[  0,\omega_{n+1}\right\rangle $. (Again, this does not involve any
induction!) The argumentation that $\boldsymbol{Y}_{n+1}$ is indeed a desired
inverse subsystem of $\boldsymbol{X}$, in the case $n+1$, is quite similar to
that of the case $n=1$. This completes the proof of the theorem.
\end{proof}

\begin{remark}
\label{R2}One has to say that the proofs of Lemma 1 and, implicitly, of
Theorem 1 are rather exceptional. Namely, it is not an often case that one
uses the ordinary and transfinite induction simultaneously such that they
alternate throughout the whole process.
\end{remark}

It remains to compare the classification of objects in $Dpro$-$\mathcal{C}$ to
that in $pro$-$\mathcal{C}$. Firstly, we need a delay-analogue of [3], Lemma I.1.2.

\begin{lemma}
\label{L2}Let $(f,f_{b})\in Dinv$-$\mathcal{C(}\boldsymbol{X},\boldsymbol{Y}%
)$, where the indexing set $B$ of $\boldsymbol{Y}$ is cofinite. Then there
exists an $(f^{\prime},f_{b}^{\prime})\in Dinv$-$\mathcal{C(}\boldsymbol{X}%
,\boldsymbol{Y})$ such that

\noindent(i) $f^{\prime}:B\rightarrow A$ is increasing;

\noindent(ii) $(f^{\prime},f_{b}^{\prime})$ is special, i,e.,

$(\forall b\in B)(\exists b_{\ast}\geq b)(\forall b^{\prime}\geq b_{\ast})$
$q_{bb^{\prime}}f_{b^{\prime}}^{\prime}=f_{b}^{\prime}p_{f^{\prime
}(b)f^{\prime}(b^{\prime})}$:

\noindent(iii) $(f^{\prime},f_{b}^{\prime})\overset{d}{\sim}(f,f_{b})$.
\end{lemma}

\begin{proof}
Since $(f,f_{b}):\emph{X}\rightarrow\emph{Y}$ belongs to $Dinv$-$\mathcal{C}$,
it follows that

$(\forall b\in B)(\exists b_{\ast}\geq b)(\forall b^{\prime}\geq b_{\ast
})(\exists a\geq f(b),f(b^{\prime}))(\forall a^{\prime}\geq a)$

$q_{bb^{\prime}}f_{b^{\prime}}p_{f(b^{\prime})a^{\prime}}=f_{b}%
p_{f(b)a^{\prime}}$.

\noindent Since $B$ is cofinite, one can define a function $u:B\rightarrow A$
by putting

$u(b)=f(b),$ $\left\vert b\right\vert =0$;

$u(b)=a$, $\left\vert b\right\vert >0$,

\noindent where $a\geq f(b_{i\ast})$, $i\in\{1,\cdots,k_{b_{\ast}}\}$,
$b_{i}\leq b_{\ast}$ (there is at most finitely predecessors of $b_{\ast}$!).
Notice that $f\leq u$. By [3], Lemma I.1.1, there exists an increasing function

$f^{\prime}:B\rightarrow A$, $u\leq f^{\prime}$.

Therefore, we have

$f_{b}p_{f(b)f^{\prime}(b)}=q_{bb^{\prime}}f_{b^{\prime}}p_{f(b^{\prime
})f^{\prime}(b^{\prime})}$

\noindent whenever $b^{\prime}\geq b_{\ast}$. Since $f\leq u\leq f^{\prime}$,
we may put

$f_{b}^{\prime}=f_{b}p_{f(b)f^{\prime}(b)}:X_{f^{\prime}(b)}\rightarrow Y_{b}%
$, $b\in B$.

\noindent It readily follows that $(f^{\prime},f_{b}^{\prime}):\boldsymbol{X}%
\rightarrow\boldsymbol{Y}$ is a delay-morphism of $Dinv$-$\mathcal{C}$
satisfying the special condition (ii). In order to verify condition (iii),
just notice that $f_{b}^{\prime}$ is a shift of $f_{b}$ (trough
$\boldsymbol{X}$).
\end{proof}

The next proposition is a delay-analogue of [3], Theorem I.1.3.

\begin{proposition}
\label{P4} For every $\boldsymbol{f}\in Dpro$-$\mathcal{C}(\boldsymbol{X}%
,\boldsymbol{Y})$ there exists an $\boldsymbol{f}^{\prime}\in Dpro$%
-$\mathcal{C}(\boldsymbol{X}^{\prime},\boldsymbol{Y}^{\prime})$ such that

\noindent(i) $\boldsymbol{X}^{\prime}$ and $\boldsymbol{Y}^{\prime}$ are
indexed over the same cofinite set

\noindent(ii) every term of $\boldsymbol{X}^{\prime}$ ($\boldsymbol{Y}%
^{\prime}$) is a term of $\boldsymbol{X}$ ($\boldsymbol{Y}$);

\noindent(iii) there exist isomorphisms $\boldsymbol{i}:\boldsymbol{X}%
\rightarrow\boldsymbol{X}^{\prime}$ and $\boldsymbol{j}:\boldsymbol{Y}%
\rightarrow\boldsymbol{Y}^{\prime}$ in $Dpro$-$\mathcal{C}$ such that
$\boldsymbol{jf}=\boldsymbol{f}^{\prime}\boldsymbol{i}$;

\noindent(iv) $\boldsymbol{f}^{\prime}$ admits a level representative
$(1_{C},f_{c}^{\prime})\in Dinv$-$\mathcal{C}(\boldsymbol{X}^{\prime
},\boldsymbol{Y}^{\prime})$.
\end{proposition}

\begin{proof}
We are following the way in the proof of [3], Theorem I.1.3. Let

$\boldsymbol{f}:\boldsymbol{X}=(X_{a},p_{aa^{\prime}},A)\rightarrow
(Y_{b},q_{bb^{\prime}},B)=\boldsymbol{Y}$

\noindent be a delay-morphism of $\boldsymbol{X}$ to $\boldsymbol{Y}$. By
Proposition 2, there is no loss of generality in assuming that $A$ and $B$ are
ordered and cofinite. By Lemma 2, there is a representative $(f,f_{b})$ of
$\boldsymbol{f}$ such that $f:B\rightarrow A$ is increasing and

$(\forall b\in B)(\exists b_{\ast}\geq b)(\forall b^{\prime}\geq b_{\ast})$
$q_{bb^{\prime}}f_{b^{\prime}}^{\prime}=f_{b}^{\prime}p_{f^{\prime
}(b)f^{\prime}(b^{\prime})}$.

\noindent Let us define a $(C,\leq)$ by putting

$c=\{c\equiv(a,b)\mid f(b)\leq a\}\subseteq A\times B$,

$(c=(a,b)\leq(a^{\prime},b^{\prime})=c^{\prime})\Leftrightarrow(a\leq
a^{\prime}\wedge b\leq b^{\prime})$.

\noindent Then $(C,\leq)$ is an ordered directed set, that is cofinite as
well, because $A$ and $B$ are cofinite. Now put

$X_{c}^{\prime}=X_{a}$, $Y_{c}^{\prime}=Y_{b}$, $p_{cc^{\prime}}^{\prime
}=p_{aa^{\prime}}$, $q_{gnc^{\prime}}^{\prime}=q_{bb^{\prime}}$.

\noindent One readily verifies that

$\boldsymbol{X}^{\prime}=(X_{c}^{\prime},p_{cc^{\prime}}^{\prime},C)$,
$\boldsymbol{Y}^{\prime}=(Y_{c}^{\prime},q_{cc^{\prime}}^{\prime},C)$

\noindent are delay-inverse systems in $\mathcal{C}.$ We now define an
$(f^{\prime},f_{c}^{\prime})$,

$f^{\prime}:C\rightarrow C$, $f_{c}^{\prime}:X_{f^{\prime}(c)}^{\prime
}\rightarrow Y_{c}^{\prime}$, $c=(a,b)\in C$,

\noindent by putting

$f^{\prime}=1_{c}$, $f_{c}^{\prime}=fl_{b}p_{f(b)a}:X_{a}\rightarrow Y_{b}$.

Let us verify that $(1_{c},f_{c}^{\prime})\in Dinv$-$\mathcal{C}%
(\boldsymbol{X}^{\prime},\boldsymbol{Y}^{\prime})$. Given a $c=(a,b)\in C$,
choose a $c_{\ast}=(a_{\ast},b_{\ast})$, where $b_{\ast}\geq b$ corresponds to
$b$ in $B$ condition on $(f,f_{b})$), while $a_{\ast}\geq f(b),f(b_{\ast
}),a^{\ast}$, where $a^{\ast}$ is a \textquotedblleft commutation
index\textquotedblright\ for $f(b_{\ast})$ (condition on $\boldsymbol{X}$).
Let $c^{\prime}=(a^{\prime},b^{\prime})\geq c_{\ast}$. Then $f_{c}^{\prime
}p_{cc^{\prime}}^{\prime}=f_{b}p_{f(b)a}p_{aa^{\prime}}=f_{b}p_{f(b)a^{\prime
}}=q_{bb^{\prime}}f_{b^{\prime}}p_{f(b^{\prime})a^{\prime}}=q_{cc^{\prime}%
}^{\prime}f_{c^{\prime}}^{\prime}$.

\noindent It remains to prove (iii). As in the proof of [3], Theorem I.1.3, put

$i:C\rightarrow A$, $i(c)=a$,

$i_{c}=1_{X_{a}}:X_{a}\rightarrow X_{a}=X_{c}^{\prime}$, $c=(a,b)\in C$,

$j:C\rightarrow B$, $j(c)=b$,

$j_{c}=1_{Y_{b}}:Y_{b}\rightarrow Y_{b}=Y_{c}^{\prime}$, $c=(a,b)\in C$.

\noindent It is almost trivial to verify that $(i,i_{c}):\boldsymbol{X}%
\rightarrow\boldsymbol{X}^{\prime}$ and $(j,j_{c}):\boldsymbol{Y}%
\rightarrow\boldsymbol{Y}^{\prime}$ are delay-morphisms of $Dinv$%
-$\mathcal{C}$. Thus $\boldsymbol{i}=[(i,i_{c})]\in Dpro$-$\mathcal{C}%
(\boldsymbol{X}^{\prime},\boldsymbol{X})$ and $\boldsymbol{j}=[(j,j_{c})]\in
Dpro$-$\mathcal{C}(\boldsymbol{Y},\boldsymbol{Y}^{\prime})$. Further, one
straightforwardly verifies that

$(1_{n},f_{c}^{\prime})(j,j_{c})\overset{d}{\sim}(j,j_{c})(f,f_{b})$,

\noindent i.e., that $\boldsymbol{f}^{\prime}\boldsymbol{i}=\boldsymbol{jf}$
holds true. Finally,a (rather tedious) verification of the statement that
$\boldsymbol{i}$ and $\boldsymbol{j}$ are isomorphisms of $Dpro$-$\mathcal{C}$
works in a way quite analogous to that in $pro$-$\mathcal{C}$
\end{proof}

We finish this note with the following two most significant facts concerning
the isomorphism classifications of the objects in $Dpro$-$\mathcal{C}$ and
$pro$-$\mathcal{C}$.

\begin{theorem}
\label{T2}Let $\boldsymbol{X}$ and $\boldsymbol{Y}$ be inverse systems in a
category $\mathcal{C}$, i.e., $\boldsymbol{X},\boldsymbol{Y}\in Ob(pro$%
-$\mathcal{C})$, having indexing sets of cardinalities $\aleph_{n}$ and
$\aleph_{n^{\prime}}$ repetitively, $n,n^{\prime}\in%
\mathbb{N}
_{0}$. If $\boldsymbol{X}\cong\boldsymbol{Y}$ in $Dpro$-$\mathcal{C}$, then
also $\boldsymbol{X}\cong\boldsymbol{Y}$ in $pro$-$\mathcal{C}$.
\end{theorem}

\begin{proof}
By Proposition 4, we may assume that $\boldsymbol{X}=(X_{c},p_{cc^{\prime}%
},C)$ and $\boldsymbol{Y}=(Y_{c},q_{cc^{\prime}},C)$, where $(C,\leq)$ is
ordered and cofinite, and that

$\boldsymbol{f}=[(1_{C},f_{c})]:\boldsymbol{X}\rightarrow\boldsymbol{Y}$

\noindent is a level isomorphism of $Dpro$-$\mathcal{C}$. In order to
construct an isomorphism of $\boldsymbol{X}$ to $\boldsymbol{Y}$ in
$pro$-$\mathbb{C}$, one has to \textquotedblleft repeat\textquotedblright\ the
whole proof of Theorem 1, now (instead for $\boldsymbol{X}\in Ob(Dinv$%
-$\mathcal{C})$) for the level representative $(1_{C},f_{c})\in Mor(Dinv$%
-$\mathcal{C})$ of $\boldsymbol{f}$. It yields a level morphism

$(1_{C^{\prime}},f_{c^{\prime}}^{\prime}=f_{c^{\prime}}):\boldsymbol{X}%
^{\prime}=(X_{c},p_{cc^{\prime}},C^{\prime})\rightarrow(Y_{c},q_{cc^{\prime}%
},C^{\prime})=\boldsymbol{Y}^{\prime}$

\noindent of $inv$-$\mathcal{C}$, where $(C^{\prime},\leq)$ is cofinal in $C$. Then

$\boldsymbol{f}^{\prime}=[(1_{C^{\prime}},f_{c})]:\boldsymbol{X}^{\prime
}\rightarrow\boldsymbol{Y}^{\prime}$

\noindent is a morphism of $pro$-$\mathcal{C}$. Since $(C^{\prime},\leq)$ is
cofinal in $(C,\leq)$ and $\boldsymbol{f}$ is an isomorphism (of
$Dpro$.$\mathcal{C})$, it readily follows that the extracted $\boldsymbol{f}%
^{\prime}$ is an isomorphism of $pro$-$\mathcal{C}$. In addition, by [3],
Theorem I.1.1, the restriction morphisms $\boldsymbol{i}:\boldsymbol{X}%
\rightarrow\boldsymbol{X}^{\prime}$ and $\boldsymbol{j}:\boldsymbol{Y}%
\rightarrow\boldsymbol{Y}^{\prime}$ are isomorphisms of $pro$-$\mathcal{C}$.
Finally, the composite

$\boldsymbol{j}^{-1}\boldsymbol{f}^{\prime}\boldsymbol{i}:\boldsymbol{X}%
\rightarrow\boldsymbol{Y}$

\noindent is an isomorphism of $pro$-$\mathcal{C}$.
\end{proof}

\begin{corollary}
\label{C1}Let $\boldsymbol{X},\boldsymbol{Y}\in Ob(Dpro$-$\mathcal{C})$ and
let $\boldsymbol{X}^{\prime},\boldsymbol{Y}^{\prime}\in Ob(pro$-$\mathcal{C})$
such that $\boldsymbol{X}\cong\boldsymbol{X}^{\prime}$ and $\boldsymbol{Y}%
\cong\boldsymbol{Y}^{\prime}$ in $Dpro$-$\mathcal{C}$. Then $\boldsymbol{X}%
\cong\boldsymbol{Y}$ in $Dpro$-$\mathcal{C}$ if and only if $\boldsymbol{X}%
^{\prime}\cong\boldsymbol{Y}^{\prime}$ in $pro$-$\mathcal{C}.$
\end{corollary}

\begin{proof}
Since $pro$-$\mathcal{C}$ is a subcategory of $Dpro$-$\mathcal{C}$, the
sufficiency part follows trivially. Conversely, since

$\boldsymbol{X}^{\prime}\boldsymbol{\cong X}\cong\boldsymbol{Y}\cong%
\boldsymbol{Y}^{\prime}$

\noindent in $Dpro-\mathcal{C}$, it follows that $\boldsymbol{X}^{\prime
}\boldsymbol{\cong Y}^{\prime}$ in $Dpro$-$\mathcal{C}$. Then, by Theorem 2,
$\boldsymbol{X}^{\prime}\boldsymbol{\cong Y}^{\prime}$ in $pro$-$\mathcal{C}$
as well.
\end{proof}

In the second part we consider the general case of an indexing set $A$ having
cardinality $\left\vert A\right\vert =\kappa$, where $\kappa>\aleph_{n}$ for
all $n\in%
\mathbb{N}
_{0}$.

AAA

\begin{center}
\textbf{References \smallskip}
\end{center}

\noindent\lbrack1] \textit{H. Herlich and G. E. Strecker,} Category Theory, An
Introduction, Allyn and Bacon Inc., Boston, 1973.

\noindent\lbrack2] \textit{S. Eilenberg and N. Steenrod,} Foundations of
Algebraic Topology, Princeton University Press, Princeton, 1952.

\noindent\lbrack3] \textit{S. Marde\v{s}i\'{c} and Jack Segal, }Shape Theory,
North Holland, Amsterdam, 1982.

\noindent\lbrack4] \textit{V. Matijevi\'{c} and L.R. Rubin,} Delay-inverse
systems, Topology Appl. \textbf{348} (2024) (10889).

\noindent\lbrack5] \textit{V. Matijevi\'{c} and L.R. Rubin,} $A$ categorical
approach to delay-inverse systems, preprint.

\noindent\lbrack6] \textit{N. Ugle\v{s}i\'{c} and B. \v{C}ervar,} The concept
of a weak shape type, International J. of Pure and Applied Math. \textbf{39}
(2007), 363-428.

\end{document}